\documentclass[twoside,11pt]{article}

\usepackage{jmlr2e}

\ShortHeadings{}{}

\usepackage{hyperref}
\usepackage{url}
\usepackage{amsmath}
\usepackage{amssymb}
\usepackage{amsfonts}
\usepackage{leftidx}

\usepackage{amsthm}
\usepackage{bm}
\usepackage[labelformat=simple]{subcaption}

\usepackage[normalem]{ulem}
\usepackage{enumitem}
\usepackage[linesnumbered]{algorithm2e}
\usepackage{multicol}

\newcommand\ci{\perp\!\!\!\perp}

\def\presuper#1#2%
  {\mathop{}%
   \mathopen{\vphantom{#2}}^{#1}%
   \kern-\scriptspace%
   #2}

\newtheorem{definition1}{Definition}
\newtheorem{theorem1}{Theorem}
\newtheorem{lemma1}{Lemma}

\usepackage{tikz}
\usepackage{tensor}
\usetikzlibrary{decorations.pathreplacing,angles,quotes}
\usepackage{tkz-tab}
\usetikzlibrary{plotmarks}

\usetikzlibrary{arrows, automata}
\usetikzlibrary{decorations.pathreplacing,angles,quotes}

\usepackage[margin=1.15in]{geometry}

\makeatletter
\newtheorem*{rep@theorem}{\rep@title}
\newcommand{\newreptheorem}[2]{%
\newenvironment{rep#1}[1]{%
 \def\rep@title{#2 \ref{##1}}%
 \begin{rep@theorem}}%
 {\end{rep@theorem}}}
\makeatother

\newreptheorem{theorem}{Theorem}
\newreptheorem{lemma}{Lemma}
\newreptheorem{proposition}{Proposition}

\makeatletter
\def\ps@pprintTitle{%
   \let\@oddhead\@empty
   \let\@evenhead\@empty
   \let\@oddfoot\@empty
   \let\@evenfoot\@oddfoot
}
\makeatother

\makeatletter
\newcommand*\rel@kern[1]{\kern#1\dimexpr\macc@kerna}
\newcommand*\widebar[1]{%
  \begingroup
  \def\mathaccent##1##2{%
    \rel@kern{0.8}%
    \overline{\rel@kern{-0.8}\macc@nucleus\rel@kern{0.2}}%
    \rel@kern{-0.2}%
  }%
  \macc@depth\@ne
  \let\math@bgroup\@empty \let\math@egroup\macc@set@skewchar
  \mathsurround\z@ \frozen@everymath{\mathgroup\macc@group\relax}%
  \macc@set@skewchar\relax
  \let\mathaccentV\macc@nested@a
  \macc@nested@a\relax111{#1}%
  \endgroup
}
\makeatother

\begin{document}

\title{The Global Markov Property for a Mixture of DAGs}

\author{\name Eric V. Strobl\\
       \addr Vanderbilt University}

\editor{TBA}

\maketitle

\begin{abstract}Real causal processes may contain feedback loops and change over time. In this paper, we model cycles and non-stationary distributions using a mixture of directed acyclic graphs (DAGs). We then study the conditional independence (CI) relations induced by a density that factorizes according to a mixture of DAGs in two steps. First, we generalize d-separation for a single DAG to \textit{mixture d-separation} for a mixture of DAGs. We then utilize the mixture d-separation criterion to derive a global Markov property that allows us to read off the CI relations induced by a mixture of DAGs using a particular summary graph. This result has potentially far reaching applications in algorithm design for causal discovery.
\end{abstract}
\begin{keywords}
Causality, Global Markov Property, Directed Acyclic Graph, Cycles
\end{keywords}

\section{The Problem}

Causal processes in nature may contain cycles and the joint density over the random variables may change over time. However, most modern representations of causality cannot accommodate cycles and non-stationarity simultaneously. For example, the \textit{directed acyclic graph} (DAG) does not contain cycles by virtue of its acyclicity \citep{Spirtes01}. Structural equation models with cycles assume a stationary distribution \citep{Spirtes94,Forre17,Strobl18}. Dynamic Bayesian networks also assume stationarity within each time step \citep{Dagum92,Dagum95}. We therefore require an alternative representation of causality in order to better model many real causal processes.

Recently, \cite{Strobl19} proposed to generalize the DAG using a mixture of DAGs. Recall that we can utilize a DAG $\mathbb{G}$ over $\bm{X}$ to represent a joint density $f(\bm{X})$ that factorizes as follows:
\begin{equation} \label{eq_DAG}
    f(\bm{X}) = \prod_{i=1}^p f(X_i | \textnormal{Pa}_{\mathbb{G}}(X_i)),
\end{equation}
where $\textnormal{Pa}_{\mathbb{G}}(X_i)$ refers to the parents, or direct causes, of $X_i \in \bm{X}$. Notice that each parent set remains fixed over time. Under the mixture of DAGs framework, we consider an auxiliary variable $T$ and assume that the joint density $f(\bm{X} \cup T)$ factorizes according to a DAG $\mathbb{G}^t$ over  $\bm{Z} = \bm{X} \cup T$ at any time point $T=t$:
\begin{equation} \nonumber
     f(\bm{X}|T=t)f(T=t)
    =\prod_{i=1}^p f(X_i | \textnormal{Pa}_{\mathbb{G}^{T=t}}(X_i),T=t) f(T=t).
\end{equation}
Notice that the parent sets may now vary with time because they are indexed by $T$. As a result, the density $f(\bm{X}|T)$ and the DAG structure over $\bm{X}$ may also change across time. We then consider a \textit{mixture of DAGs} by ``mixing across time'' or integrating out $T$:
\begin{equation} \label{eq_fac_mix}
    f(\bm{X}) = \sum_T f(\bm{X}|T)f(T)
    =\sum_T \Big[ \prod_{i=1}^p f(X_i | \textnormal{Pa}_{\mathbb{G}^T}(X_i),T)f(T) \Big].
\end{equation}

The above mixture of DAGs representation allows us to model both cycles and non-stationarity simultaneously. Consider for example the causal process shown in Figure 1 involving two random variables. Intuitively, a cycle occurs when we iteratively ``cycle through'' or ``unravel'' the variables in the feedback loop. In this case, first $X_i$ causes $X_j$, then $X_j$ causes $X_i$, then $X_i$ causes $X_j$ and so forth. This suggests that we can decompose the cycle into two DAGs $X_i \rightarrow X_j$ and $X_j \rightarrow X_i$. Suppose for simplicity that $X_i$ causes $X_j$ at time point 1, and $X_j$ causes $X_i$ at time point 2. The joint density $f(X_i,X_j)$ therefore factorizes as follows at $T=1$ and $T=2$, respectively:
\begin{equation} \nonumber
\begin{aligned}
&f(X_i,X_j|T=1) = f(X_j|X_i, T = 1)f(X_i| T = 1),\\
&f(X_i,X_j|T=2) = f(X_i|X_j, T = 2)f(X_j| T = 2).\\
\end{aligned}
\end{equation}
Notice then that $f(\bm{X}|T) \not = f(\bm{X})$; we thus say that $f(\bm{X})$ is \textit{non-stationary}. Unfortunately, we may not always have the luxury of sampling from only one of the above two densities. We therefore instead suppose that we can at least sample from both time steps, or from the following mixture density: 
\begin{equation} \nonumber
f(X_i,X_j|T=1)f(T=1) + f(X_i,X_j|T=2)f(T=2).
\end{equation}
We have included a dataset sampled according to the above mixture density in Figure 1. There, the  two red samples correspond to samples from $f(X_i,X_j|T=1)$ and the three black samples correspond to samples from $f(X_i,X_j|T=2)$. We now have a clear understanding of what it means to sample from the cycle in Figure 1; we have samples obtained from a mixture of two densities that each factorize according to a DAG because sometimes we obtain samples when $X_i$ causes $X_j$, and other times we obtain samples when $X_j$ causes $X_i$.

\begin{figure}
\centering
\begin{subfigure}{0.155\textwidth}
\centering
\resizebox{\linewidth}{!}{
\begin{tikzpicture}[scale=1.0, shorten >=1pt,auto,node distance=2.8cm, semithick]
                    
\tikzset{vertex/.style = {inner sep=0.4pt}}
\tikzset{edge/.style = {->,> = latex'}}
 
\node[vertex] (1) at  (0,0) {$X_i$};
\node[vertex] (2) at  (1.5,0) {$X_j$};

\draw[edge, bend right] (1) to (2);
\draw[edge, bend right] (2) to (1);
\end{tikzpicture}
}
\caption{}  \label{fig_cycle}
\end{subfigure}
\begin{subtable}{0.45\textwidth}
  \centering
\begin{tabular}{ c c}
  $X_i$ & $X_j$ \\
	\hline
  \color{red}{0.21} & \color{red}{-0.20} \\
  0.68 & -0.47 \\
  \color{red}{1.05} & \color{red}{-0.19} \\
  0.72 & -1.40 \\
  0.13 & -0.56 \\ 	
  \hline
\end{tabular}
\caption{} \label{table_SB1}
\end{subtable}
\caption{We consider the cyclic causal process depicted in (a). We decompose (a) into two DAGs $X_i \rightarrow X_j$ and $X_j \rightarrow X_i$ at time points 1 and 2, respectively. We draw the red samples in the table in (b) from $f(X_i,X_j|T=1)$ which factorizes according to $X_i \rightarrow X_j$. Similarly, we draw the black samples in the table from $f(X_i,X_j|T=2)$ which factorizes according to $X_j \rightarrow X_i$.}
\end{figure}
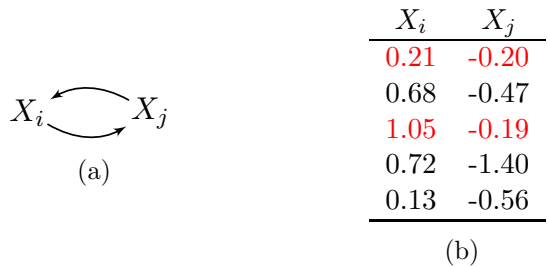

In this report, we focus on deriving the \textit{global Markov property} for a mixture of DAGs. This property will allow us to read off the conditional independence (CI) relations implied by the joint density in Equation \eqref{eq_fac_mix} using a graphical criterion. Note that two attempts have been made to derive this property in the past. \cite{Spirtes94} for example derived a global Markov property for a mixture of DAGs, but his proposal misses many important CI relations lying on directed paths. \cite{Strobl19} also made an attempt to derive the property in order to account for the additional CI relations, but several investigators discovered a counter-example to his proposal (see Acknowledgements). The above failures thus highlight the difficulty of deriving the property in the general setting. However, the past attempts suggested a proof strategy, which we utilized herein to successfully derive the global Markov property from first principles. 

\section{Preliminaries}

We now provide additional background knowledge to keep this report self-contained.

\subsection{Graphical Terminology}

We will consider both undirected and directed graphs. An \textit{undirected graph} represents each variable in $\bm{X}$ as a vertex and contains \textit{undirected edges} ``$-$'' between the vertices. We say that two vertices $X_i$ and $X_j$ are \textit{adjacent} in an undirected graph when we have $X_i - X_j$. Consider three disjoint subsets of $\bm{X}$, denoted as $\bm{A}$, $\bm{B}$ and $\bm{C}$. A set of vertices $\bm{A}$ forms a \textit{clique} in an undirected graph when all of the vertices in $\bm{A}$ are adjacent to each other. We call a sequence of undirected edges between $\bm{A}$ and $\bm{B}$ that does not pass through a vertex more than once as an \textit{undirected path} between $\bm{A}$ and $\bm{B}$. We say that $\bm{A}$ and $\bm{B}$ are \textit{connected} given $\bm{C}$ if and only if there exists an undirected path between $\bm{A}$ and $\bm{B}$ that does not pass through $\bm{C}$. Moreover, $\bm{A}$ and $\bm{B}$ are \textit{separated} given $\bm{C}$ if and only if $\bm{A}$ and $\bm{B}$ are not connected given $\bm{C}$. An undirected path $\Pi$ is \textit{active} between $\bm{A}$ and $\bm{B}$ given $\bm{C}$ when $\Pi$ lies between $\bm{A}$ and $\bm{B}$ but does not pass through $\bm{C}$. We say that a joint density over $\bm{X}$ satisfies the \textit{global Markov property} with respect to (w.r.t.) an undirected graph when the following property holds: if $\bm{A}$ and $\bm{B}$ are separated given $\bm{C}$, then $\bm{A}$ and $\bm{B}$ are conditionally independent given $\bm{C}$, denoted as $\bm{A} \ci \bm{B} | \bm{C}$ for shorthand.

We will summarize a causal process using directed graphs. A \textit{directed graph} $\mathbb{G}$ represents each variable in $\bm{X}$ as a vertex and contains \textit{directed edges} ``$\rightarrow$'' or ``$\leftarrow$'' between the vertices. We say that $X_i$ is a \textit{direct cause} or \textit{parent} of $X_j$, when we have the directed edge $X_i \rightarrow X_j$. We denote this relation as $X_i \in \textnormal{Pa}_{\mathbb{G}}(X_j)$ for shorthand. Similarly, $X_j$ is a \textit{child} of $X_i$. Two vertices $X_i$ and $X_j$ are \textit{adjacent}, if there exists a directed edge between $X_i$ and $X_j$ irrespective of its direction. We write $X_i * \!\! - \!\! * X_j$ when we have either $X_i \rightarrow X_j$ or $X_i \leftarrow X_j$. If we have $X_i \rightarrow X_j \leftarrow X_k$ in a directed graph, then we refer to $X_j$ as a \textit{collider}. We also call $X_i$ and $X_k$ \textit{spouses}. On the other hand, the vertex $X_j$ is a \textit{non-collider}, if we have $X_i \rightarrow X_j \rightarrow X_k$, $X_i \leftarrow X_j \leftarrow X_k$ or $X_i \leftarrow X_j \rightarrow X_k$. We call a sequence of directed edges between $\bm{A}$ and $\bm{B}$ that does not pass through a vertex more than once as a \textit{directed path} between $\bm{A}$ and $\bm{B}$. Next, $X_i$ is an \textit{ancestor} of $X_j$, denoted as $X_i \in \textnormal{Anc}_{\mathbb{G}}(X_j)$, when there exists a directed path from $X_i$ to $X_j$. We also apply the definition of an ancestor to a set of vertices $\bm{A}$ as follows: $\textnormal{Anc}(\bm{A}) = \{X_i | X_i \in \textnormal{Anc}(X_j) \text{ for some } X_j \in \bm{A}\}$. A directed graph contains a \textit{cycle} or a \textit{feedback loop} when $X_i$ is an ancestor of $X_j$, and we have the directed edge $X_j \rightarrow X_i$. A directed graph is more specifically called a \textit{directed acyclic graph} (DAG), if the directed graph does not contain cycles. We say that $\bm{A}$ and $\bm{B}$ are \textit{d-connected} given $\bm{C}$ if and only if there exists a directed path $\Pi$ between $\bm{A}$ and $\bm{B}$ such that the following two criteria hold: (a) $X_k \in \bm{W}$ for every collider $X_k$ on $\Pi$, and (b) $X_k \not \in \bm{W}$ for every non-collider $X_k$ on $\Pi$. The vertices $\bm{A}$ and $\bm{B}$ are \textit{d-separated} given $\bm{C}$ if and only if $\bm{A}$ and $\bm{B}$ are not d-connected given $\bm{C}$. We denote d-connection and d-separation as $\bm{A} \not \ci_d \bm{B} | \bm{C}$ and $\bm{A} \ci_d \bm{B} | \bm{C}$, respectively. The following lemma allows us to combine d-connecting paths to form new d-connecting paths:
\begin{lemma1} \label{lem_d_conn} (Lemma 2.5 in \cite{Colombo12})
Suppose that $X_i$ and $X_j$ are not in $\bm{W} \subseteq \bm{X} \setminus \{X_i, X_j\}$, there is a sequence $\sigma$ of distinct vertices in $\bm{X}$ from $X_i$ to $X_j$, and there is a set $\mathcal{P}$ of paths such that:
\begin{enumerate}
\item for each pair of adjacent vertices $X_v$ and $X_w$ in $\sigma$, there is a unique path in $\mathcal{P}$ that d-connects $X_v$ and $X_w$ given $\bm{W}$; 
\item if a vertex $X_q$ in $\sigma$ is in $\bm{W}$, then the paths in $\mathcal{P}$ that contain $X_q$ as an endpoint collide at $X_q$;
\item if for three vertices $X_v$, $X_w$ and $X_q$ occurring in that order in $\sigma$, the d-connecting paths in $\mathcal{P}$ between $X_v$ and $X_w$, and between $X_w$ and $X_q$ collide at $X_w$, then $X_w \in \textnormal{Anc}(\bm{W})$.
\end{enumerate}
Then there is a path $\Pi_{X_i X_j}$ in $\mathbb{G}$ that d-connects $X_i$ and $X_j$ given $\bm{W}$.
\end{lemma1}
\noindent We say that a joint density over $\bm{X}$ satisfies the \textit{global Markov property} w.r.t. a directed graph when the following property holds: if $\bm{A} \ci_d \bm{B} | \bm{C}$, then $\bm{A} \ci \bm{B} | \bm{C}$ \citep{Lauritzen90}. We can associate a \textit{moral graph} to a directed graph by \textit{marrying} the spouses, or drawing an undirected edge between the spouses, and then converting all directed edges into undirected edges. We will utilize the following equivalence relation between d-separation in a directed graph and connection in the corresponding moral graph:
\begin{lemma1} \label{lem_moral}
(Proposition 5.13 on page 72 in \cite{Cowell99}) $\bm{A} \ci_d \bm{B} | \bm{C}$ if and only if $\bm{A}$ and $\bm{B}$ are connected given $\bm{C}$ in the moral graph over $\textnormal{Anc}(\bm{A} \cup \bm{B} \cup \bm{C})$.
\end{lemma1}

\subsection{Further Details on Mixture of DAGs}

Recall that we can utilize a DAG to represent a joint density that factorizes according to Equation 
\eqref{eq_DAG}. We can also generalize the DAG to a mixture of DAGs in order to model the joint density in Equation \eqref{eq_fac_mix}. Unfortunately, the joint density in Equation \eqref{eq_fac_mix} does not imply any CI relations over $\bm{X}$ because all of the variables are children of $T$.

The set $\textnormal{Pa}_{\mathbb{G}^T}(X_i)$ may however not vary on the support of $f(T)$. We may also have $f(X_i | \textnormal{Pa}_{\mathbb{G}^T}(X_i)\setminus T,T) = f(X_i | \textnormal{Pa}_{\mathbb{G}^T}(X_i)\setminus T)$. In other words, the parent set remains unchanged and the conditional density does not vary across time. Let $\bm{X}^\emptyset \subseteq \bm{X}$ denote the set of variables satisfying the above two criteria. We can then write the following for those $X_i \in \bm{X}^\emptyset$:
\begin{align}
f(X_i|\textnormal{Pa}_{\mathbb{G}^T}(X_i),T) =& f(X_i|\textnormal{Pa}_{\mathbb{G}^T}(X_i)) \nonumber\\
=&f(X_i|\textnormal{Pa}_{\mathbb{G}^\emptyset}(X_i)), \nonumber
\end{align}
where $\textnormal{Pa}_{\mathbb{G}^\emptyset}(X_i)$ denotes a parent set that does not vary over time. We may now rewrite Equation \eqref{eq_fac_mix} as follows:
\begin{align}
&\sum_T \Big[f(T)\prod_{i=1}^p f(X_i|\textnormal{Pa}_{\mathbb{G}^T}(X_i),T)\Big] \nonumber\\
=&\sum_T \Big[f(T)\prod_{i=1}^r f(X_i|\textnormal{Pa}_{\mathbb{G}^T}(X_i))\prod_{i=1}^u f(X_i|\textnormal{Pa}_{\mathbb{G}^\emptyset}(X_i))\Big]. \label{eq_fac_endT}
\end{align}
where we assume that $T \in \textnormal{Pa}_{\mathbb{G}^T}(X_i)$ for all $X_i \in [\bm{X} \setminus \bm{X}^\emptyset]$, and $T \not \in \textnormal{Pa}_{\mathbb{G}^\emptyset}(X_i)$ for all $X_i \in \bm{X}^\emptyset$.

\subsection{Mother Graph}

Let $\mathcal{G}$ refer to the set of unique DAGs over $\bm{X} \cup T$ indexed by $T$. Note that $\mathbb{G}^\emptyset \in \mathcal{G}$. \cite{Strobl19} introduced the notion of a mother graph $\mathcal{M}$ to graphically represent a mixture of DAGs. The mother graph is a DAG formed by plotting all of the DAGs in $\mathcal{G}$ next to each other. We therefore refer to each graph in $\mathcal{G}$ as a \textit{sub-DAG} of $\mathcal{M}$. Note that $\mathcal{G}$ can only contain a finite number of DAGs because the number of possible DAGs over a finite random vector $\bm{Z}$ is finite. Let $q \in \mathbb{N}^+$ denote the number of DAGs in $\mathcal{G}$. We use the notation $\mathcal{G}_i$ to refer to the $i^{\textnormal{th}}$ sub-DAG of $\mathcal{M}$.

Suppose for instance that $\mathcal{G}$ contains three sub-DAGs $X_i \rightarrow X_j \leftarrow X_k$, $X_i \leftarrow X_j \hspace{4mm} X_k$ and $X_i \hspace{4mm} X_j \leftarrow X_k$ with $T$ a direct cause of $X_j$. We can then plot these three DAGs next to each other as in Figure \ref{fig_mother} to create $\mathcal{M}$. The superscripts in $\mathcal{M}$ index the DAGs in $\mathcal{G}$. Let $\bm{D},\bm{E}$ and $\bm{F}$ denote disjoint subsets of $\bm{Z}$. Also let $Z_i^\prime = \{Z_i^1,\dots, Z_i^q\}$ denote the set of \textit{vertices} in $\mathcal{M}$ corresponding to the \textit{variable} $Z_i$. Similarly, let $\bm{D}^\prime = \cup_{Z_i \in \bm{D}} Z_i^\prime$. We say that the two sets of variables $\bm{D}$ and $\bm{E}$ are \textit{d-separated} given $\bm{F}$ in $\mathcal{M}$ if and only if the vertices $\bm{D}^\prime$ and $\bm{E}^\prime$ are d-separated given $\bm{F}^\prime$ in $\mathcal{M}$. We for instance have $X_i \ci_d X_k$ in $\mathcal{M}$ in Figure \ref{fig_mother} because we have $\{X_i^1,X_i^2,X_i^3\} \ci_d \{X_k^1,X_k^2,X_k^3\}$. From here on, we refer to the vertices $X_i^\prime$ in $\mathcal{M}$ using the variable $X_i$ and drop the superscripts in the graphical representation of $\mathcal{M}$ in order to simplify notation. Similarly, we refer to the set of vertices $\bm{D}^\prime$ with the set of variables $\bm{D}$. We write $Z_i *\!\! - \!\!* Z_j$ in $\mathcal{M}$ if and only if we have $Z_i^k *\!\! - \!\!* Z_j^k$ in $\mathcal{M}$ for some $\mathcal{G}_k \in \mathcal{G}$. We say that $Z_i$ is a \textit{parent} of $Z_j$ in $\mathcal{M}$ if and only if $Z_i^k$ is a parent of $Z_j^k$ in $\mathcal{M}$ for some $\mathcal{G}_k \in \mathcal{G}$. The same holds for \textit{children}, \textit{spouses} and \textit{ancestors}.

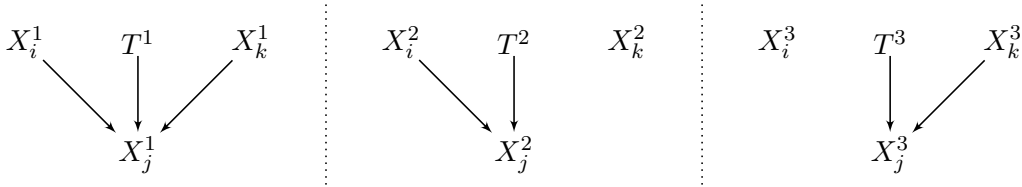
\begin{figure} \label{fig_mother}
\centering{
\begin{tikzpicture}[scale=1.0, shorten >=1pt,auto,node distance=2.8cm, semithick]
                    
\tikzset{vertex/.style = {inner sep=0.4pt}}
\tikzset{edge/.style = {->,> = latex'}}
 
\node[vertex] (1) at  (0,0) {$X_i^1$};
\node[vertex] (2) at  (1.5,-1.5) {$X_j^1$};
\node[vertex] (3) at  (3,0) {$X_k^1$};
\node[vertex] (4) at  (1.5,0) {$T^1$};

\draw[edge] (1) to (2);
\draw[edge] (4) to (2);
\draw[edge] (3) to (2);

\draw[dotted] (4,0.5) to (4,-2);

\node[vertex] (5) at  (5,0) {$X_i^2$};
\node[vertex] (6) at  (6.5,-1.5) {$X_j^2$};
\node[vertex] (7) at  (8,0) {$X_k^2$};
\node[vertex] (8) at  (6.5,0) {$T^2$};

\draw[edge] (5) to (6);
\draw[edge] (8) to (6);

\draw[dotted] (9,0.5) to (9,-2);

\node[vertex] (5) at  (10,0) {$X_i^3$};
\node[vertex] (6) at  (11.5,-1.5) {$X_j^3$};
\node[vertex] (7) at  (13,0) {$X_k^3$};
\node[vertex] (8) at  (11.5,0) {$T^3$};

\draw[edge] (7) to (6);
\draw[edge] (8) to (6);
\end{tikzpicture}
}
\caption{An example of a mother graph $\mathcal{M}$. The superscripts correspond to $\mathcal{G}_1, \mathcal{G}_2$ and $\mathcal{G}_3$.} \label{fig_mother}
\end{figure}

\section{New Definitions}

\begin{figure}
\centering
\begin{subfigure}{0.47\textwidth}
\centering
\begin{tikzpicture}[scale=1.0, shorten >=1pt,auto,node distance=2.8cm, semithick]
                    
\tikzset{vertex/.style = {inner sep=0.4pt}}
\tikzset{edge/.style = {->,> = latex'}}
 
\node[vertex] (1) at  (0,0) {$X_i$};
\node[vertex] (2) at  (1.5,-1.5) {$X_j$};
\node[vertex] (3) at  (3,0) {$X_k$};
\node[vertex] (4) at  (1.5,0) {$T$};

\draw[edge] (1) to (2);
\draw[edge] (4) to (2);

\draw[dotted] (4,0.5) to (4,-2);

\node[vertex] (5) at  (5,0) {$X_i$};
\node[vertex] (6) at  (6.5,-1.5) {$X_j$};
\node[vertex] (7) at  (8,0) {$X_k$};
\node[vertex] (8) at  (6.5,0) {$T$};

\draw[edge] (8) to (6);
\draw[edge] (7) to (6);
\end{tikzpicture}
 \caption{} \label{fig_CE}
\end{subfigure}
\vspace{10mm}

\begin{subfigure}{0.47\textwidth}
\centering{
\begin{tikzpicture}[scale=1.0, shorten >=1pt,auto,node distance=2.8cm, semithick]
                    
\tikzset{vertex/.style = {inner sep=0.4pt}}
\tikzset{edge/.style = {->,> = latex'}}
 
\node[vertex] (1) at  (0,0) {$X_i$};
\node[vertex] (2) at  (1.5,-1.5) {$X_j$};
\node[vertex] (3) at  (3,0) {$X_k$};
\node[vertex] (4) at  (1.5,0) {$T$};

\draw[edge] (2) to (1);
\draw[edge] (4) to (2);
\draw[edge] (4) to (1);

\draw[dotted] (4,0.5) to (4,-2);

\node[vertex] (5) at  (5,0) {$X_i$};
\node[vertex] (6) at  (6.5,-1.5) {$X_j$};
\node[vertex] (7) at  (8,0) {$X_k$};
\node[vertex] (8) at  (6.5,0) {$T$};

\draw[edge] (8) to (6);
\draw[edge] (7) to (6);
\draw[edge] (8) to (5);
\end{tikzpicture}
} \caption{} \label{fig_md_sep}
\end{subfigure}
\caption{Subfigure (a) provides an example of both an m-collider $X_j$ and an m-path between $X_i$ and $X_k$. We also have examples of mother graphs where (a) $X_i \not \ci_{md} X_k | X_j$ and (b) $X_i \ci_{md} X_k$. Notice that d-separation finds an erroneous CI relation in (a), whereas m-d-separation does not.}
\end{figure}
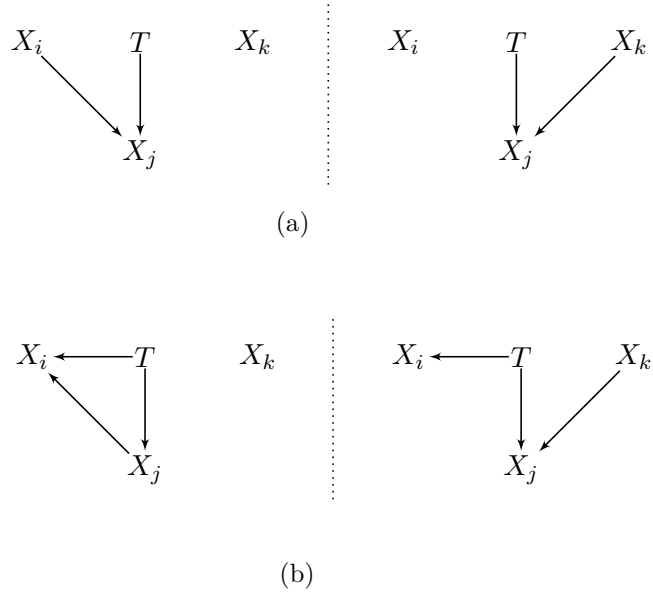

We first require some new definitions about the mother graph before we can state the main result. We reserve $Z_{p+1}=T$ so that $Z_i = X_i$ for any $i \leq p$. We define the \textit{m-collider}, a generalization of a collider in a single directed graph:
\begin{definition1} \label{def_mcollider}
The variable $X_j$ is a \textit{mixture collider} (m-collider) in $\mathcal{M}$ if and only if at least one of the following conditions holds for the triple $\langle Z_i, X_j, Z_k \rangle$:
\begin{enumerate}
    \item $Z_i \rightarrow X_j \leftarrow  Z_k$ in any sub-DAG of $\mathcal{M}$;
    \item $X_i \rightarrow X_j \leftarrow  T$ in one sub-DAG of $\mathcal{M}$ and $T \rightarrow X_j \leftarrow  X_k$ in another sub-DAG.
\end{enumerate}
\end{definition1}
\noindent Notice that the first part of the definition corresponds to the definition of a collider in a single directed graph. We require the second condition in order to account for conditional dependence relations that may be induced between sub-DAGs in $\mathcal{M}$. Figures \ref{fig_mother} and \ref{fig_CE} provide examples of m-colliders.

We also require the definition of an \textit{m-path} in a mother graph. Let $\bm{D},\bm{E}$ and $\bm{F}$ correspond to disjoint subsets of $\bm{Z}$.
\begin{definition1}
A \textit{mixture path} (m-path) $\Pi$ exists between $\bm{D}$ and $\bm{E}$ in $\mathcal{M}$ if and only if there exists a sequence of triples between $\bm{D}$ and $\bm{E}$ such that at least one of the following two conditions holds for each triple $\langle Z_i,Z_j,Z_k \rangle$:
\begin{enumerate}
    \item $Z_i * \!\! - \!\!* Z_j * \!\! - \!\! * Z_k$ exists in at least one sub-DAG of $\mathcal{M}$;
    \item $X_i  \rightarrow X_j \leftarrow  T$ in one sub-DAG of $\mathcal{M}$ and $T  \rightarrow X_j \leftarrow X_k$ in another sub-DAG.
\end{enumerate}
\end{definition1}
\noindent Notice that the second condition in the above definition corresponds to the second condition in the definition of an m-collider. There thus exists an m-path between $X_i$ and $X_k$ in the mother graph shown in Figure \ref{fig_CE}, even though there does not exist a directed path between $X_i$ and $X_k$ in any sub-DAG of $\mathcal{M}$.

We are now ready to define \textit{m-d-connection} as follows:
\begin{definition1} \label{def_md}
We say that $\bm{D}$ and $\bm{E}$ are m-d-connected given $\bm{F}$, denoted as $\bm{D} \not \ci_{md} \bm{E} | \bm{F}$, if and only if there exists an m-path $\Pi$ between $\bm{D}$ and $\bm{E}$ such that the following two conditions hold:
\begin{enumerate}
    \item $X_k \in \bm{F}$ for every m-collider $X_k$ on $\Pi$;
    \item $Z_k \not \in \bm{F}$ for every non-m-collider $Z_k$ on $\Pi$.
\end{enumerate}
\end{definition1}
\noindent Notice then that we have $X_i \not \ci_{md} X_j | X_k$ in both Figure \ref{fig_CE} and Figure \ref{fig_mother}.  We say that $\bm{D}$ and $\bm{E}$ are m-d-separated given $\bm{F}$, denoted as $\bm{D} \ci_{md} \bm{E} | \bm{F}$, if and only if they are not m-d-connected given $\bm{F}$. Figure \ref{fig_md_sep} provides an example of m-d-separation because we have $X_i \ci_{md} X_k$. We can relate d-separation and m-d-separation via the following lemma:
\begin{lemma1} \label{lem_md_d}
If $\bm{D} \ci_{md} \bm{E} | \bm{F}$ in $\mathcal{M}$, then $\bm{D} \ci_{d} \bm{E} | \bm{F}$ in $\mathcal{M}$.
\end{lemma1}
\begin{proof}
We prove this by contrapositive. If $\bm{D} \not \ci_{d} \bm{E} | \bm{F}$ in $\mathcal{M}$, then the following conditions hold for at least one directed path $\Pi$ between $\bm{D}$ and $\bm{E}$:
\begin{enumerate}
    \item $X_k \in \bm{F}$ for every collider $X_k$ on $\Pi$;
    \item $Z_k \not \in \bm{F}$ for every non-collider $X_k$ on $\Pi$.
\end{enumerate}
Notice that a directed path $\Pi$ between $\bm{D}$ and $\bm{E}$ must also be an m-path between $\bm{D}$ and $\bm{E}$ by definition. It follows that the following two conditions hold for at least one m-path $\Pi$ between $\bm{D}$ and $\bm{E}$:
\begin{enumerate}
    \item $X_k \in \bm{F}$ for every m-collider $X_k$ on $\Pi$;
    \item $Z_k \not \in \bm{F}$ for every non-m-collider $Z_k$ on $\Pi$ because, if $X_k$ is an m-collider on $\Pi$ and the second condition in Definition \ref{def_mcollider} holds but the first condition does not, then $\Pi$ is not a directed path between $\bm{D}$ and $\bm{E}$.
\end{enumerate}
 We therefore conclude that we have $\bm{D} \not \ci_{md} \bm{E} | \bm{F}$ in $\mathcal{M}$.
\end{proof}

\section{Main Result}

We are now ready to derive the global Markov property for a mixture of DAGs. Recall that $\bm{A},\bm{B}$ and $\bm{C}$ denote disjoint subsets of $\bm{X}$. We have:
\begin{theorem1} \label{thm_GMP}
(Global Markov Property) If $\bm{A} \ci_{md} \bm{B} | \bm{C}$ in $\mathcal{M}$, then $\bm{A} \ci \bm{B} | \bm{C}$.
\end{theorem1}
\begin{proof}
We first consider $\widebar{\mathcal{M}}$, the moral graph of $\mathcal{M}$. Let $\widebar{\mathcal{M}}_{\textnormal{Anc}_{\mathcal{M}}(\bm{A}\cup\bm{B} \cup \bm{C} )}$ denote the moral graph of the ancestral set containing $\bm{A}\cup\bm{B} \cup \bm{C}$. We then consider a partition of variables (not just vertices) $\ddot{\bm{A}} \cup \ddot{\bm{B}} \cup \bm{C} = \textnormal{Anc}_{\mathcal{M}}(\bm{A}\cup\bm{B} \cup \bm{C} )$ such that $\bm{A} \subseteq \ddot{\bm{A}}$, $\bm{B} \subseteq \ddot{\bm{B}}$, and $\ddot{\bm{A}}$, $\ddot{\bm{B}}$ and $\bm{C}$ are disjoint sets of variables. We also require that $\ddot{\bm{A}}$ and $\ddot{\bm{B}}$ be separated by $\bm{C}$ in $\widebar{\mathcal{M}}_{\textnormal{Anc}_{\mathcal{M}}(\bm{A}\cup\bm{B} \cup \bm{C} )}$. We now argue that such a partition is possible. By Lemma \ref{lem_md_d}, we have $\bm{A} \ci_{d} \bm{B} | \bm{C}$ in $\mathcal{M}$, so $\bm{A}$ and  $\bm{B}$ are also separated by $\bm{C}$ in $\widebar{\mathcal{M}}_{\textnormal{Anc}_{\mathcal{M}}(\bm{A}\cup\bm{B} \cup \bm{C} )}$ by Lemma \ref{lem_moral}. Now consider the set of variables $\bm{H} = \textnormal{Anc}_{\mathcal{M}}(\bm{A}\cup\bm{B} \cup \bm{C} )\setminus (\bm{A} \cup \bm{B} \cup \bm{C}).$ We have two situations for each $Z_i \in \bm{H}$ in $\widebar{\mathcal{M}}_{\textnormal{Anc}_{\mathcal{M}}(\bm{A}\cup\bm{B} \cup \bm{C} )}$:
\begin{enumerate}
    \item There does not exist an undirected path between $Z_i$ and $\bm{A}$ for every sub-DAG or an undirected path between $Z_i$ and $\bm{B}$ for every sub-DAG (or both) that is active given $\bm{C}$. More specifically:
    \begin{enumerate}
    \item If there does not exist an undirected path between $Z_i$ and $\bm{A}$ that is active given $\bm{C}$ for every sub-DAG, but such a path does exist between $Z_i$ and $\bm{B}$ for some sub-DAG, then place $Z_i$ in $\ddot{\bm{B}}$.
    \item If there does not exist an undirected path between $Z_i$ and $\bm{B}$ that is active given $\bm{C}$ for every sub-DAG, but such a path does exist between $Z_i$ and $\bm{A}$ for some sub-DAG, then place $Z_i$ in $\ddot{\bm{A}}$.
    \item If there does not exist an undirected path between $Z_i$ and $\bm{A}$ that is active given $\bm{C}$ for every sub-DAG and there likewise does not exist such a path between $Z_i$ and $\bm{B}$ for every sub-DAG, then place $Z_i$ in either $\ddot{\bm{A}}$ or $\ddot{\bm{B}}$ (but not both).
    \end{enumerate}
    \item There exists an undirected path between $Z_i$ and $\bm{A}$ for some sub-DAG and an undirected path between $Z_i$ and $\bm{B}$ for some sub-DAG that are both active given $\bm{C}$. If the two sub-graphs in $\widebar{\mathcal{M}}_{\textnormal{Anc}_{\mathcal{M}}(\bm{A}\cup\bm{B} \cup \bm{C} )}$ correspond to the same sub-DAG in $\mathcal{M}$, then $\bm{A}$ and  $\bm{B}$ would be connected given $\bm{C}$ in $\widebar{\mathcal{M}}_{\textnormal{Anc}_{\mathcal{M}}(\bm{A}\cup\bm{B} \cup \bm{C} )}$; a contradiction. Now suppose that the two sub-graphs in $\widebar{\mathcal{M}}_{\textnormal{Anc}_{\mathcal{M}}(\bm{A}\cup\bm{B} \cup \bm{C} )}$ corresponds to different sub-DAGs in $\mathcal{M}$, denoted as $\mathcal{G}_1$ and $\mathcal{G}_2$. Note that $Z_i$ and $\bm{A}$ are d-connected given $\bm{C}$ in $\mathcal{G}_1$, so there exists a d-connecting path $\Pi_{Z_i\bm{A}}$ between $Z_i$ and $\bm{A}$ given $\bm{C}$ in $\mathcal{G}_1$. As a result, every vertex on $\Pi_{Z_i\bm{A}}$ is a member of $\textnormal{Anc}_{\mathcal{M}}(Z_i \cup \bm{A} \cup \bm{C} )$. Similarly, $Z_i$ and $\bm{B}$ are d-connected given $\bm{C}$ in $\mathcal{G}_2$, so there exists a d-connecting path $\Pi_{Z_i\bm{B}}$ between $Z_i$ and $\bm{B}$ given $\bm{C}$ in $\mathcal{G}_2$. Every vertex on $\Pi_{Z_i\bm{B}}$ is a member of $\textnormal{Anc}_{\mathcal{M}}(Z_i \cup \bm{B} \cup \bm{C} )$. Because $\Pi_{Z_i\bm{A}}$ exists in $\mathcal{G}_1$ but not in $\mathcal{G}_2$ and likewise $\Pi_{Z_i\bm{B}}$ exists in $\mathcal{G}_2$ but not in $\mathcal{G}_1$, there must exist a child of $T$ on $\Pi_{Z_i\bm{A}}$ and a child of $T$ on $\Pi_{Z_i\bm{B}}$, so $T \in \textnormal{Anc}_{\mathcal{M}}(Z_i\cup\bm{A} \cup \bm{B} \cup \bm{C} )$. Since $Z_i \in \textnormal{Anc}_{\mathcal{M}}(\bm{A} \cup \bm{B} \cup \bm{C} )$, it follows that we may more specifically claim $T \in \textnormal{Anc}_{\mathcal{M}}(\bm{A} \cup \bm{B} \cup \bm{C} )$. Let $\Pi_{T\bm{A}}$ denote a shortest d-connecting path between $T$ and $\bm{A}$ on $\Pi_{Z_i\bm{A}}$, and likewise let $\Pi_{T\bm{B}}$ denote a shortest d-connecting path between $T$ and $\bm{B}$ on $\Pi_{Z_i\bm{B}}$. Setting $\mathcal{P} = \{\Pi_{T\bm{A}}, \Pi_{T\bm{B}}\} $ and then invoking Lemma \ref{lem_d_conn}, we form a d-connecting path between $\bm{A}$ and $\bm{B}$ given $\bm{C}$, which implies that $\bm{A}$ and $\bm{B}$ are connected given $\bm{C}$ in $\widebar{\mathcal{M}}_{\textnormal{Anc}_{\mathcal{M}}(\bm{A}\cup\bm{B} \cup \bm{C} )}$; another contradiction. We can therefore conclude that there does not exist an undirected path between $Z_i$ and $\bm{A}$ for some sub-DAG and an undirected path between $Z_i$ and $\bm{B}$ for some sub-DAG that are both active given $\bm{C}$, even if the two sub-graphs correspond to different sub-DAGs.
\end{enumerate}
We conclude that there exists a disjoint partition of variables  $\ddot{\bm{A}} \cup \ddot{\bm{B}} \cup \bm{C} = \textnormal{Anc}_{\mathcal{M}}(\bm{A}\cup\bm{B} \cup \bm{C} ) $. Moreover, by the impossibility of the second point above, $\ddot{\bm{A}}$ and $\ddot{\bm{B}}$ must be separated by $\bm{C}$ in $\widebar{\mathcal{M}}_{\textnormal{Anc}_{\mathcal{M}}(\bm{A}\cup\bm{B} \cup \bm{C} )}$.

Now assume that $T \not \in \textnormal{Anc}_{\mathcal{M}}(\bm{A}\cup\bm{B} \cup \bm{C} )$. We may then consider all of the cliques in $\widebar{\mathcal{M}}_{\textnormal{Anc}_{\mathcal{M}}(\bm{A}\cup\bm{B} \cup \bm{C} )}$ corresponding to each vertex and its married parents. Denote this set of cliques as $\mathcal{E}$. Also let $\mathcal{E}_{\ddot{\bm{B}}}$ denote the set of cliques in $\mathcal{E}$ that have non-empty intersection with $\ddot{\bm{B}}$. By the above paragraph, $\ddot{\bm{A}}$ and $\ddot{\bm{B}}$ are non-adjacent in $\widebar{\mathcal{M}}_{\textnormal{Anc}_{\mathcal{M}}(\bm{A}\cup\bm{B} \cup \bm{C} )}$, so no clique in $\mathcal{E}_{\ddot{\bm{B}}}$ can contain a member of $\ddot{\bm{A}}$. We also have $\ddot{\bm{B}} \cap e =\emptyset$ for all $e \in \mathcal{E} \setminus \mathcal{E}_{\ddot{\bm{B}}}$. We can write the following using $\widebar{\mathcal{M}}_{\textnormal{Anc}_{\mathcal{M}}(\bm{A}\cup\bm{B} \cup \bm{C} )}$:
\begin{equation} \nonumber
\begin{aligned}
    f(\ddot{\bm{A}}, \ddot{\bm{B}}, \bm{C})&= \prod_{\{\bm{D}_i \cup \textnormal{Pa}_{\mathcal{M}}(\bm{D}_i)\} \in \mathcal{E} \setminus \mathcal{E}_{\ddot{\bm{B}}}} f(\bm{D}_i | \textnormal{Pa}_{\mathcal{M}}(\bm{D}_i) ) \prod_{\{\bm{E}_i \cup \textnormal{Pa}_{\mathcal{M}}(\bm{E}_i)\} \in  \mathcal{E}_{\ddot{\bm{B}}}} f(\bm{E}_i | \textnormal{Pa}_{\mathcal{M}}(\bm{E}_i) ) \\
    &= \prod_{e \in \mathcal{E} \setminus \mathcal{E}_{\ddot{\bm{B}}}} \gamma(e) \prod_{e \in \mathcal{E}_{\ddot{\bm{B}}}} \gamma(e) = \gamma(\ddot{\bm{A}}, \bm{C} )  \gamma(\ddot{\bm{B}}, \bm{C} ),
\end{aligned}
\end{equation}
where $\gamma$ denotes a non-negative function. We then proceed by integrating out $[\ddot{\bm{A}}\cup\ddot{\bm{B}}] \setminus [\bm{A}\cup\bm{B}]$:
\begin{equation} \nonumber
\begin{aligned}
f(\bm{A}, \bm{B}, \bm{C}) &=\sum_{[\ddot{\bm{A}}\cup\ddot{\bm{B}}] \setminus [\bm{A}\cup\bm{B}]} f(\ddot{\bm{A}}, \ddot{\bm{B}}, \bm{C}) = \sum_{[\ddot{\bm{A}} \setminus \bm{A}] \cup [\ddot{\bm{B}} \setminus \bm{B}]} f(\ddot{\bm{A}}, \ddot{\bm{B}}, \bm{C}) \\
&= \sum_{[\ddot{\bm{A}} \setminus \bm{A}] \cup [\ddot{\bm{B}} \setminus \bm{B}]} \gamma(\ddot{\bm{A}}, \bm{C} ) \gamma(\ddot{\bm{B}}, \bm{C}) = \Big[\sum_{[\ddot{\bm{B}} \setminus \bm{B}]} \Big[ \sum_{[\ddot{\bm{A}} \setminus \bm{A}]} \gamma(\ddot{\bm{A}}, \bm{C}) \Big] \gamma(\ddot{\bm{B}}, \bm{C})\Big]\\
&= \sum_{[\ddot{\bm{A}} \setminus \bm{A}]} \gamma(\ddot{\bm{A}}, \bm{C})  \sum_{[\ddot{\bm{B}} \setminus \bm{B}]} \gamma(\ddot{\bm{B}}, \bm{C}) = \gamma(\bm{A}, \bm{C}) \gamma(\bm{B}, \bm{C}),
\end{aligned}
\end{equation}
where the fifth equality follows because $[\ddot{\bm{A}} \setminus \bm{A}] \cap [\ddot{\bm{B}} \setminus \bm{B}]=\emptyset$ by construction. The conclusion follows by the sixth equality in this case.

Now assume that $T \in \textnormal{Anc}_{\mathcal{M}}(\bm{A}\cup\bm{B} \cup \bm{C} )$. Since $\bm{A}, \bm{B}$ and $\bm{C}$ are disjoint subsets of $\bm{X}$, it follows that $T$ is contained in either $\ddot{\bm{A}} \setminus \bm{A}$ or $\ddot{\bm{B}} \setminus \bm{B}$ (but not both). Assume without loss of generality that $T$ is contained in $\ddot{\bm{B}} \setminus \bm{B}$ and therefore $\ddot{\bm{B}}$. As a result, $\ddot{\bm{A}}$ and therefore $\bm{A}$ cannot contain a child of $T$ because this would imply that $\ddot{\bm{A}}$ and $\ddot{\bm{B}}$ are connected given $\bm{C}$ via $T$ in $\widebar{\mathcal{M}}_{\textnormal{Anc}_{\mathcal{M}}(\bm{A}\cup\bm{B} \cup \bm{C} )}$. It follows that only $\ddot{\bm{B}}$ or $\bm{C}$ (or both) can contain the children of $T$. As a side note, notice that connection given $\bm{C}$  in $\widebar{\mathcal{M}}_{\textnormal{Anc}_{\mathcal{M}}(\bm{A}\cup\bm{B} \cup \bm{C} )}$ alone (or equivalently d-separation via Lemma \ref{lem_moral}) does not exclude the possibility that $\bm{A}$ and hence $\ddot{\bm{A}}$ contains a spouse of $T$ (e.g., Figure \ref{fig_CE} where $\bm{A}=X_i,\bm{B}=X_k$, $\bm{C}=X_j$, $\ddot{\bm{B}} = X_k \cup T$ and $X_i \ci_d X_k | X_j$). We will nevertheless see that m-d-separation gracefully handles this situation as detailed below. 

Note that $\ddot{\bm{A}}$ and $ \ddot{\bm{B}}$ are non-adjacent $\widebar{\mathcal{M}}_{\textnormal{Anc}_{\mathcal{M}}(\bm{A}\cup\bm{B} \cup \bm{C} )}$ because they are separated given $\bm{C}$. We then have two scenarios:
\begin{enumerate}
    \item (Part 1) $\ddot{\bm{B}}$ contains at least one child of $T$. Then, $\ddot{\bm{A}}$ cannot contain a spouse of $T$ because this would imply that $\ddot{\bm{A}}$ and $\ddot{\bm{B}}$ are adjacent in $\widebar{\mathcal{M}}_{\textnormal{Anc}_{\mathcal{M}}(\bm{A}\cup\bm{B} \cup \bm{C} )}$. Hence, $\ddot{\bm{A}}$ and $T$ are non-adjacent in $\widebar{\mathcal{M}}_{\textnormal{Anc}_{\mathcal{M}}(\bm{A}\cup\bm{B} \cup \bm{C} )}$, and we may more specifically claim that $\ddot{\bm{A}}$ and $ \ddot{\bm{B}} \cup T$ are non-adjacent in $\widebar{\mathcal{M}}_{\textnormal{Anc}_{\mathcal{M}}(\bm{A}\cup\bm{B} \cup \bm{C} )}$.
    
    Next let $\mathcal{E}_{\ddot{\bm{B}} \cup T}$ denote the set of cliques in $\mathcal{E}$ that have non-empty intersection with $\ddot{\bm{B}} \cup T$. By the above paragraph, $\ddot{\bm{A}}$ and $\ddot{\bm{B}}\cup T$ are non-adjacent in $\widebar{\mathcal{M}}_{\textnormal{Anc}_{\mathcal{M}}(\bm{A}\cup\bm{B} \cup \bm{C} )}$, so no clique in $\mathcal{E}_{\ddot{\bm{B}} \cup T}$ can contain a member of $\ddot{\bm{A}}$. We also have $(\ddot{\bm{B}} \cup T) \cap e =\emptyset$ for all $e \in \mathcal{E} \setminus \mathcal{E}_{\ddot{\bm{B}} \cup T}$. 

Note that $\mathcal{E}_{\ddot{\bm{B}} \cup T}$ contains all of the children and spouses of $T$ (that are both in $\ddot{\bm{A}} \cup \ddot{\bm{B}}$ and not in $\ddot{\bm{A}} \cup \ddot{\bm{B}}$, i.e. in $\bm{C}$). We next take $\bm{C}$ and consider the partition $\bm{C} \cup [\bm{C} \setminus \widehat{\bm{C}}]$, where $\widehat{\bm{C}}$ corresponds to the vertices in both $\bm{C}$ and a clique in $\mathcal{E}_{\ddot{\bm{B}} \cup T}$. Next, we consider the density $f(\ddot{\bm{A}}, \ddot{\bm{B}}, \bm{C})$. We can write the following using $\widebar{\mathcal{M}}_{\textnormal{Anc}_{\mathcal{M}}(\bm{A}\cup\bm{B} \cup \bm{C} )}$:
\begin{equation} \nonumber
\begin{aligned}
    &\hspace{5mm}f(\ddot{\bm{A}}, \ddot{\bm{B}}, \bm{C})= \sum_T f(\ddot{\bm{A}}, \ddot{\bm{B}}, \bm{C}, T)\\ 
    &= \prod_{\{\bm{D}_i \cup \textnormal{Pa}_{\mathcal{M}}(\bm{D}_i)\} \in \mathcal{E} \setminus \mathcal{E}_{\ddot{\bm{B}} \cup T}} f(\bm{D}_i | \textnormal{Pa}_{\mathcal{M}}(\bm{D}_i) ) \sum_T\prod_{\{\bm{E}_i \cup \textnormal{Pa}_{\mathcal{M}}(\bm{E}_i)\} \in  \mathcal{E}_{\ddot{\bm{B}} \cup T }} f(\bm{E}_i | \textnormal{Pa}_{\mathcal{M}}(\bm{E}_i) ) \\
    &= \prod_{e \in \mathcal{E} \setminus \mathcal{E}_{\ddot{\bm{B}} \cup T}} \gamma(e) \sum_T \prod_{e \in \mathcal{E}_{\ddot{\bm{B}} \cup T}} \gamma(e) = \gamma(\ddot{\bm{A}}, [\bm{C} \setminus \widehat{\bm{C}}] ) \sum_T \gamma(\ddot{\bm{B}}, \widehat{\bm{C}}, T )\\
    &= \gamma(\ddot{\bm{A}}, [\bm{C} \setminus \widehat{\bm{C}}] ) \gamma(\ddot{\bm{B}}, \widehat{\bm{C}} ),
\end{aligned}
\end{equation}
where $\gamma$ denotes a non-negative function. The summation after the second equality is possible because $\mathcal{E}_{\ddot{\bm{B}} \cup T}$ contains every clique that intersects with $T$ by construction.

We then proceed by integrating out $[\ddot{\bm{A}}\cup\ddot{\bm{B}}] \setminus [\bm{A}\cup\bm{B}] = [\ddot{\bm{A}} \setminus \bm{A}] \cup [\ddot{\bm{B}} \setminus \bm{B}]$:
\begin{equation} \nonumber
\begin{aligned}
f(\bm{A}, \bm{B}, \bm{C})
&= \sum_{[\ddot{\bm{A}} \setminus \bm{A}] \cup [\ddot{\bm{B}} \setminus \bm{B}]} f(\ddot{\bm{A}}, \ddot{\bm{B}}, \bm{C}) = \sum_{[\ddot{\bm{A}} \setminus \bm{A}] \cup [\ddot{\bm{B}} \setminus \bm{B}]} \gamma(\ddot{\bm{A}}, [\bm{C} \setminus \widehat{\bm{C}}] ) \gamma(\ddot{\bm{B}}, \widehat{\bm{C}}) \\
&= \sum_{[\ddot{\bm{A}} \setminus \bm{A}]} \gamma(\ddot{\bm{A}}, [\bm{C} \setminus \widehat{\bm{C}}])  \sum_{[\ddot{\bm{B}} \setminus \bm{B}]} \gamma(\ddot{\bm{B}}, \widehat{\bm{C}}) = \gamma(\bm{A}, [\bm{C} \setminus \widehat{\bm{C}}]) \gamma(\bm{B}, \widehat{\bm{C}})\\
&= \gamma(\bm{A}, \bm{C}) \gamma(\bm{B}, \bm{C}).
\end{aligned}
\end{equation}
The conclusion follows by the fifth equality in this case.

    \item $\ddot{\bm{B}}$ does not contain any children of $T$. Then all of the children of $T$ in $\textnormal{Anc}_{\mathcal{M}}(\bm{A}\cup\bm{B} \cup \bm{C} )$ must lie in $\bm{C}$. We have four sub-cases: 
    \begin{enumerate}
        \item $\ddot{\bm{B}}$ contains a spouse of $T$ in $\textnormal{Anc}_{\mathcal{M}}(\bm{A}\cup\bm{B} \cup \bm{C} )$ but $\ddot{\bm{A}}$ does not. Since $\ddot{\bm{A}}$ does not contain a spouse of $T$ or a child of $T$, we know that $\ddot{\bm{B}} \cup T$ and $\ddot{\bm{A}}$ are non-adjacent in $\widebar{\mathcal{M}}_{\textnormal{Anc}_{\mathcal{M}}(\bm{A}\cup\bm{B} \cup \bm{C} )}$. It follows that we can proceed just like Part 1 using $\mathcal{E}_{\ddot{\bm{B}} \cup T}$ and $\mathcal{E} \setminus \mathcal{E}_{\ddot{\bm{B}} \cup T}$. 
        
        \item $\ddot{\bm{A}}$ contains a spouse of $T$ in $\textnormal{Anc}_{\mathcal{M}}(\bm{A}\cup\bm{B} \cup \bm{C} )$ but $\ddot{\bm{B}}$ does not. We proceed similarly to the argument of (a). Since $\ddot{\bm{B}}$ does not contain a spouse of $T$ or a child of $T$, we know that $\ddot{\bm{A}} \cup T$ and $\ddot{\bm{B}}$ are non-adjacent in $\widebar{\mathcal{M}}_{\textnormal{Anc}_{\mathcal{M}}(\bm{A}\cup\bm{B} \cup \bm{C} )}$. It follows that we can proceed just like Part 1, except with the modification that we consider the cliques $\mathcal{E}_{\ddot{\bm{A}} \cup T}$ and $\mathcal{E} \setminus \mathcal{E}_{\ddot{\bm{A}} \cup T}$ as opposed to $\mathcal{E}_{\ddot{\bm{B}} \cup T}$ and $\mathcal{E} \setminus \mathcal{E}_{\ddot{\bm{B}} \cup T}$. 
        
        \item Both $\ddot{\bm{A}}$ and $\ddot{\bm{B}}$ contain spouses of $T$ in $\textnormal{Anc}_{\mathcal{M}}(\bm{A}\cup\bm{B} \cup \bm{C} )$. This implies that $\ddot{\bm{A}} \not \ci_{md} \ddot{\bm{B}} | \bm{C}$. We already know that $\bm{A} \ci_{md} \bm{B} | \bm{C}$, so $\bm{A}$ and $\bm{B}$ cannot both contain spouses of $T$. We may thus claim that one of the following three situations holds:
        \begin{enumerate}
            \item $\bm{A}$ does not contain a spouse of $T$ but $\bm{B}$ does, so that $\ddot{\bm{A}} \setminus \bm{A}$ and $\bm{B}$ contain spouses of $T$. Note further that $\bm{A}$ cannot be adjacent to $\ddot{\bm{A}} \setminus \bm{A}$ in any sub-DAG where $\ddot{\bm{A}} \setminus \bm{A}$ is a spouse of $T$ because, if this did happen, then we would have $\bm{A} \not \ci_{md} \bm{B} | \bm{C}$ in $\mathcal{M}$. It follows that $\bm{A}$ and $(\ddot{\bm{A}} \setminus \bm{A}) \cup \ddot{\bm{B}} \cup T$ are non-adjacent $\widebar{\mathcal{M}}_{\textnormal{Anc}_{\mathcal{M}}(\bm{A}\cup\bm{B} \cup \bm{C} )}$. We can therefore consider the cliques $\mathcal{E}_{(\ddot{\bm{A}} \setminus \bm{A}) \cup \ddot{\bm{B}} \cup T}$ as well as $\mathcal{E} \setminus \mathcal{E}_{(\ddot{\bm{A}} \setminus \bm{A}) \cup \ddot{\bm{B}} \cup T}$, where the latter contains all of the variables in $\bm{A}$, and proceed as in Part 1.
            \item $\bm{B}$ does not contain a spouse of $T$ but $\bm{A}$ does, so that $\bm{A}$ and $\ddot{\bm{B}} \setminus \bm{B}$ contain spouses of $T$. The argument for this proceeds analogously as the above argument. We therefore consider the cliques $\mathcal{E}_{ \ddot{\bm{A}} \cup (\ddot{\bm{B}} \setminus \bm{B}) \cup T}$ as well as $\mathcal{E} \setminus \mathcal{E}_{ \ddot{\bm{A}} \cup (\ddot{\bm{B}} \setminus \bm{B}) \cup T}$, where the latter contains all of the variables in $\bm{B}$, and proceed as in Part 1.
            \item Both $\bm{A}$ and $\bm{B}$ do not contain spouses of $T$, so that $\ddot{\bm{A}} \setminus \bm{A}$ and $\ddot{\bm{B}} \setminus \bm{B}$ contain spouses of $T$. It follows that both $\bm{A}$ and $\bm{B}$ are non-adjacent to $T$ in $\widebar{\mathcal{M}}_{\textnormal{Anc}_{\mathcal{M}}(\bm{A}\cup\bm{B} \cup \bm{C} )}$ (and non-adjacent to each other). This implies that the cliques in $\mathcal{E}_{\bm{A}}$ have empty intersection with the cliques in $\mathcal{E}_{\bm{B}}$. Let $\ddot{\bm{E}}$ correspond to those variables in $\ddot{\bm{A}} \setminus \bm{A}$ that have non-empty intersection with the cliques in $\mathcal{E}_{\bm{A}}$. Similarly let $\ddot{\bm{F}}$ correspond to those variables in $\ddot{\bm{B}} \setminus \bm{B}$ that have non-empty intersection with the cliques in $\mathcal{E}_{\bm{B}}$. Finally, let $\ddot{\bm{G}}$ correspond to those variables in $[\ddot{\bm{A}} \setminus \bm{A}] \cup [\ddot{\bm{B}} \setminus \bm{B}]$ that have non-empty intersection with the cliques in $\mathcal{E} \setminus [\mathcal{E}_{\bm{A}} \cup \mathcal{E}_{\bm{B}}]$. Notice that $\ddot{\bm{E}}, \ddot{\bm{F}}$ and $\ddot{\bm{G}}$ are disjoint subsets of $[\ddot{\bm{A}} \setminus \bm{A}] \cup [\ddot{\bm{B}} \setminus \bm{B}]$. We then write:
            \begin{equation} \nonumber
            \begin{aligned}
                f(\ddot{\bm{A}}, \ddot{\bm{B}}, \bm{C}) &= \sum_T f(\ddot{\bm{A}}, \ddot{\bm{B}}, \bm{C}, T)\\ 
                &= \prod_{\{\bm{A}_i \cup \textnormal{Pa}_{\mathcal{M}}(\bm{A}_i)\} \in \mathcal{E}_{\bm{A}}} f(\bm{A}_i | \textnormal{Pa}_{\mathcal{M}}(\bm{A}_i) ) \prod_{\{\bm{B}_i \cup \textnormal{Pa}_{\mathcal{M}}(\bm{B}_i)\} \in \mathcal{E}_{\bm{B}}} f(\bm{B}_i | \textnormal{Pa}_{\mathcal{M}}(\bm{B}_i) ) \\&\hspace{10mm}\sum_T\prod_{\{\bm{D}_i \cup \textnormal{Pa}_{\mathcal{M}}(\bm{D}_i)\} \in  \mathcal{E} \setminus [\mathcal{E}_{\bm{A}} \cup \mathcal{E}_{\bm{B}}]} f(\bm{D}_i | \textnormal{Pa}_{\mathcal{M}}(\bm{D}_i) ) \\
                &= \prod_{e \in \mathcal{E}_{\bm{A}}} \gamma(e) \prod_{e \in \mathcal{E}_{\bm{B}}} \gamma(e) \sum_T \prod_{e \in \mathcal{E} \setminus [\mathcal{E}_{\bm{A}} \cup \mathcal{E}_{\bm{B}}]} \gamma(e)\\
                &= \gamma(\bm{A}, \ddot{\bm{E}},[\bm{C} \setminus \widehat{\bm{C}}] )\gamma(\bm{B},\ddot{\bm{F}}, [\bm{C} \setminus \widehat{\bm{C}}] ) \sum_T \gamma(\ddot{\bm{G}}, \widehat{\bm{C}}, T ),
            \end{aligned}
            \end{equation}
            where $\widehat{\bm{C}}$ now corresponds to those vertices in $\widebar{\mathcal{M}}_{\textnormal{Anc}_{\mathcal{M}}(\bm{A}\cup\bm{B} \cup \bm{C} )}$ in both $\bm{C}$ and $\mathcal{E} \setminus [\mathcal{E}_{\bm{A}} \cup \mathcal{E}_{\bm{B}}]$. Integrating out $[\ddot{\bm{A}}\cup\ddot{\bm{B}}] \setminus [\bm{A}\cup\bm{B}] = [\ddot{\bm{A}} \setminus \bm{A}] \cup [\ddot{\bm{B}} \setminus \bm{B}]$, we obtain:
            \begin{equation} \nonumber
            \begin{aligned}
            f(\bm{A}, \bm{B}, \bm{C}) &= \sum_{[\ddot{\bm{A}} \setminus \bm{A}] \cup [\ddot{\bm{B}} \setminus \bm{B}]} f(\ddot{\bm{A}}, \ddot{\bm{B}}, \bm{C}) \\
            &= \sum_{[\ddot{\bm{A}} \setminus \bm{A}] \cup [\ddot{\bm{B}} \setminus \bm{B}]} \gamma(\bm{A}, \ddot{\bm{E}},[\bm{C} \setminus \widehat{\bm{C}}] )\gamma(\bm{B},\ddot{\bm{F}}, [\bm{C} \setminus \widehat{\bm{C}}] ) \sum_T \gamma(\ddot{\bm{G}}, \widehat{\bm{C}}, T ) \\
            &= \sum_{\ddot{\bm{E}}} \gamma(\bm{A}, \ddot{\bm{E}},[\bm{C} \setminus \widehat{\bm{C}}] )\Big[\sum_{\ddot{\bm{F}}} \gamma(\bm{B},\ddot{\bm{F}}, [\bm{C} \setminus \widehat{\bm{C}}] ) \sum_{\ddot{\bm{G}}\cup T} \gamma(\ddot{\bm{G}}, \widehat{\bm{C}}, T ) \Big]\\
            &= \gamma(\bm{A}, [\bm{C} \setminus \widehat{\bm{C}}]) \gamma(\bm{B}, \bm{C})= \gamma(\bm{A}, \bm{C}) \gamma(\bm{B}, \bm{C}),
            \end{aligned}
            \end{equation}
            where the third equality follows because $\ddot{\bm{E}}, \ddot{\bm{F}}$ and $\ddot{\bm{G}}$ are disjoint subsets of $[\ddot{\bm{A}} \setminus \bm{A}] \cup [\ddot{\bm{B}} \setminus \bm{B}]$ by construction.

\end{enumerate}
      \item Both $\ddot{\bm{A}}$ and $\ddot{\bm{B}}$ do not contain spouses of $T$ in $\textnormal{Anc}_{\mathcal{M}}(\bm{A}\cup\bm{B} \cup \bm{C} )$. This implies that $T$ is not adjacent to $\ddot{\bm{A}} \cup \ddot{\bm{B}}$ in  $\widebar{\mathcal{M}}_{\textnormal{Anc}_{\mathcal{M}}(\bm{A}\cup\bm{B} \cup \bm{C} )}$. We can therefore write:
                 \begin{equation} \nonumber
            \begin{aligned}
                f(\ddot{\bm{A}}, \ddot{\bm{B}}, \bm{C}) &= \sum_T f(\ddot{\bm{A}}, \ddot{\bm{B}}, \bm{C}, T)\\ 
                &= \prod_{\{\bm{D}_i \cup \textnormal{Pa}_{\mathcal{M}}(\bm{D}_i)\} \in \mathcal{E}_{\ddot{\bm{A}}}} f(\bm{D}_i | \textnormal{Pa}_{\mathcal{M}}(\bm{D}_i) ) \prod_{\{\bm{E}_i \cup \textnormal{Pa}_{\mathcal{M}}(\bm{E}_i)\} \in \mathcal{E}_{\ddot{\bm{B}}}} f(\bm{E}_i | \textnormal{Pa}_{\mathcal{M}}(\bm{E}_i) ) \\&\hspace{10mm}\sum_T\prod_{\{\bm{F}_i \cup \textnormal{Pa}_{\mathcal{M}}(\bm{F}_i)\} \in  \mathcal{E} \setminus [\mathcal{E}_{\ddot{\bm{A}}} \cup \mathcal{E}_{\ddot{\bm{B}}}]} f(\bm{F}_i | \textnormal{Pa}_{\mathcal{M}}(\bm{F}_i) ) \\
                &= \prod_{e \in \mathcal{E}_{\ddot{\bm{A}}}} \gamma(e) \prod_{e \in \mathcal{E}_{\ddot{\bm{B}}}} \gamma(e) \sum_T \prod_{e \in \mathcal{E} \setminus [\mathcal{E}_{\ddot{\bm{A}}} \cup \mathcal{E}_{\ddot{\bm{B}}}]} \gamma(e)\\
                &= \gamma(\ddot{\bm{A}}, [\bm{C} \setminus \widehat{\bm{C}}] )\Big[\gamma(\ddot{\bm{B}}, [\bm{C} \setminus \widehat{\bm{C}}] ) \sum_T \gamma( \widehat{\bm{C}}, T )\Big]\\
                &= \gamma(\ddot{\bm{A}}, \bm{C} )\gamma(\ddot{\bm{B}}, \bm{C}),
            \end{aligned}
            \end{equation}
            where $\widehat{\bm{C}}$ now corresponds to those vertices in $\widebar{\mathcal{M}}_{\textnormal{Anc}_{\mathcal{M}}(\bm{A}\cup\bm{B} \cup \bm{C} )}$ in both $\bm{C}$ and in a clique in $\mathcal{E} \setminus [\mathcal{E}_{\ddot{\bm{A}}} \cup \mathcal{E}_{\ddot{\bm{B}}}]$. Integrating out $[\ddot{\bm{A}}\cup\ddot{\bm{B}}] \setminus [\bm{A}\cup\bm{B}]=[\ddot{\bm{A}} \setminus \bm{A}] \cup [\ddot{\bm{B}} \setminus \bm{B}]$, we obtain:
            \begin{equation} \nonumber
            \begin{aligned}
             f(\bm{A}, \bm{B}, \bm{C})
            &= \sum_{[\ddot{\bm{A}} \setminus \bm{A}] \cup [\ddot{\bm{B}} \setminus \bm{B}]} f(\ddot{\bm{A}}, \ddot{\bm{B}}, \bm{C}) = \sum_{[\ddot{\bm{A}} \setminus \bm{A}] \cup [\ddot{\bm{B}} \setminus \bm{B}]} \gamma(\ddot{\bm{A}}, \bm{C} )\gamma(\ddot{\bm{B}}, \bm{C})\\ &= \sum_{\ddot{\bm{A}} \setminus \bm{A} } \gamma(\ddot{\bm{A}}, \bm{C} )\sum_{\ddot{\bm{B}} \setminus \bm{B}}\gamma(\ddot{\bm{B}}, \bm{C}) = \gamma(\bm{A}, \bm{C}) \gamma(\bm{B}, \bm{C}).
            \end{aligned}
            \end{equation}
            
\end{enumerate}
\end{enumerate}

\noindent The conclusion follows because we have considered all possible cases. 
\end{proof}
\noindent Notice that the proof of the above theorem proceeds by taking the moral graph of an ancestral set of $\mathcal{M}$, and then parsing this moral graph into cliques for all possible cases. The m-d-separation criterion always allows us to separate the cliques into two or three groups with one containing $T$ and the others not. The proof then proceeds by marginalization, allowing us to factorize the joint density $f(\bm{A},\bm{B},\bm{C})$ into a product of two non-negative functions involving $\bm{A} \cup \bm{C}$ and $\bm{B} \cup \bm{C}$.

Note that standard d-separation fails in the proof at 2 (c) because both $\bm{A}$ and $\bm{B}$ can contain spouses of $T$. Consider for example the mother graph in Figure \ref{fig_CE}. Let $\mathcal{T}_1$ denote the time points where the joint density factorizes according to $\mathcal{G}_1$ and $\mathcal{T}_2$ to $\mathcal{G}_2$. Here, $X_i \not \ci X_k | X_j$ because we can write the joint density as:
\begin{equation} \nonumber
\begin{aligned}
    f(X_i,X_j,X_k) &= \sum_{T \in \mathcal{T}_1} f(X_i)f(X_j|X_i,T)f(X_k)f(T)\\ &\hspace{5mm}+ \sum_{T \in \mathcal{T}_2}f(X_i)f(X_j|X_k,T)f(X_k)f(T),
\end{aligned}
\end{equation}
which does \textit{not} factorize into two non-negative functions $\gamma(X_i,X_j) \gamma(X_k,X_j)$. However, notice that $X_i \ci_d X_k |X_j$ in $\mathcal{M}$, so that d-separation implies an erroneous CI relation. Now let $\ddot{\bm{A}}=\bm{A}=X_i$, $\bm{B}=X_k, \ddot{\bm{B}}= X_k \cup T$ and $\bm{C} = X_j$. We can trace the problem of d-separation within the proof of Theorem \ref{thm_GMP} back to the fact that both $\ddot{\bm{A}}$ and $\ddot{\bm{B}}$ contain spouses of $T$ as in 2 (c), but so do both $\bm{A}$ and $\bm{B}$. As a result, we cannot factorize the density into cliques and ultimately into two non-negative functions involving $\bm{A} \cup \bm{C}$ and $\bm{B} \cup \bm{C}$. In contrast, the m-d-separation criterion prevents $\bm{A}$ and $\bm{B}$ from simultaneously containing spouses of $T$ and therefore leads to the correct conclusion. We can check this claim in the above example because we have $X_i \not\ci_{md} X_k |X_j$.

\section{Examples}

We now illustrate the global Markov property w.r.t. a mixture of DAGs using a few examples. We first consider the mother graph in Figure \ref{fig_md_sep}, where $X_i \ci_{md} X_k$. We may write:
\begin{equation} \nonumber
\begin{aligned}
    f(X_i,X_k) &= \sum_{X_j} \Big[ \sum_{T \in \mathcal{T}_1} f(X_k)f(X_i|X_j,T)f(X_j|T)f(T) + \sum_{T \in \mathcal{T}_2} f(X_i|T)f(X_k)f(T) \Big]\\
              &= f(X_k)\sum_{X_j} \Big[\sum_{T \in \mathcal{T}_1}f(X_i|X_j,T)f(X_j|T)f(T) + \sum_{T \in \mathcal{T}_2}f(X_i|T)f(T)\Big]\\
              &= \gamma(X_k)\gamma(X_i).
\end{aligned}  
\end{equation}
We thus conclude that m-d-separation successfully detects the independence relation $X_i \ci X_k$ in this case. 

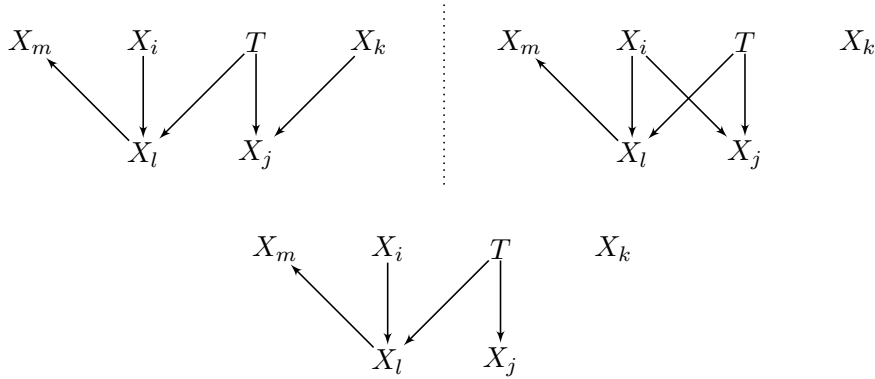
\begin{figure}
\centering{
\begin{tikzpicture}[scale=1.0, shorten >=1pt,auto,node distance=2.8cm, semithick]
                    
\tikzset{vertex/.style = {inner sep=0.4pt}}
\tikzset{edge/.style = {->,> = latex'}}
 
\node[vertex] (1) at  (0,0) {$X_i$};
\node[vertex] (5) at  (0,-1.5) {$X_l$};
\node[vertex] (6) at  (-1.5,0) {$X_m$};
\node[vertex] (2) at  (1.5,-1.5) {$X_j$};
\node[vertex] (3) at  (3,0) {$X_k$};
\node[vertex] (4) at  (1.5,0) {$T$};

\draw[edge] (4) to (2);
\draw[edge] (4) to (5);
\draw[edge] (3) to (2);
\draw[edge] (1) to (5);
\draw[edge] (5) to (6);

\draw[dotted] (4,0.5) to (4,-2);

\node[vertex] (5) at  (5,0) {$X_m$};
\node[vertex] (4) at  (6.5,-1.5) {$X_l$};
\node[vertex] (6) at  (8,0) {$T$};
\node[vertex] (1) at  (6.5,0) {$X_i$};
\node[vertex] (2) at  (8,-1.5) {$X_j$};
\node[vertex] (3) at  (9.5,0) {$X_k$};

\draw[edge] (1) to (2);
\draw[edge] (1) to (4);
\draw[edge] (6) to (2);
\draw[edge] (6) to (4);
\draw[edge] (4) to (5);

\node[vertex] (5) at  (1.75,-2.75) {$X_m$};
\node[vertex] (4) at  (3.25,-4.25) {$X_l$};
\node[vertex] (6) at  (4.75,-2.75) {$T$};
\node[vertex] (1) at  (3.25,-2.75) {$X_i$};
\node[vertex] (2) at  (4.75,-4.25) {$X_j$};
\node[vertex] (3) at  (6.25,-2.75) {$X_k$};

\draw[edge] (1) to (4);
\draw[edge] (6) to (2);
\draw[edge] (6) to (4);
\draw[edge] (4) to (5);

\end{tikzpicture}
}
\caption{A more complicated example illustrating the global Markov property, where we have $X_m \ci_{md} X_k | \{X_j, X_l\}$.} \label{fig_ex_compl}
\end{figure}

Figure \ref{fig_ex_compl} provides a more complicated example involving both an m-collider and a non-m-collider. Here, we have three sub-DAGs in $\mathcal{M}$. Notice that $X_m \ci_{md} X_k | \{X_j, X_l\}$. We may also write:
\begin{equation} \nonumber
\begin{aligned}
        f(X_j,X_k,X_l,X_m) &= \sum_{X_i} \Big[ \sum_{T \in \mathcal{T}_1} f(X_m|X_l)f(X_l|X_i,T)f(X_i)f(T)f(X_j|X_k,T)f(X_k)\\ 
            &\hspace{15mm}+ \sum_{T \in \mathcal{T}_2} f(X_m|X_l)f(X_l|X_i,T)f(X_i)f(T)f(X_j|X_i,T)f(X_k)\\ 
            &\hspace{15mm} + \sum_{T \in \mathcal{T}_3} f(X_m|X_l)f(X_l|X_i,T)f(X_i)f(T)f(X_j|T)f(X_k)\Big]\\
              &= f(X_m|X_l) \sum_{X_i} \Big[\sum_{T \in \mathcal{T}_1} f(X_l|X_i,T)f(X_i)f(T)f(X_j|X_k,T)f(X_k)\\ 
              &\hspace{15mm}+ \sum_{T \in \mathcal{T}_2} f(X_l|X_i,T)f(X_i)f(T)f(X_j|X_i,T)f(X_k)\\ 
              &\hspace{15mm} + \sum_{T \in \mathcal{T}_3} f(X_l|X_i,T)f(X_i)f(T)f(X_j|T)f(X_k)\Big]\\
              &= \gamma(X_m,X_j,X_l)\gamma(X_k,X_j,X_l),
\end{aligned}  
\end{equation}
so that we have $X_m \ci X_k | \{X_j, X_l\}$ as expected from the global Markov property.

\section{Conclusion}

We derived the global Markov property for a mixture of the DAGs. Our derivation builds upon the ideas introduced in \citep{Spirtes94} and \citep{Strobl19}. We introduced the notion of m-d-separation which in turn utilizes the notion of an m-collider. The m-collider allows us to prove the property by factorizing the joint density into two non-negative functions, one involving $T$ and the other not. We ultimately hope that this global Markov property will help investigators design algorithms for causal discovery using a mixture of DAGs in order to handle feedback loops and non-stationarity simultaneously. 

\section{Acknowledgements}
We thank Basil Saeed, Snigdha Panigrahi and Caroline Uhler for discovering a counter-example to Theorem \ref{thm_GMP} in \citep{Strobl19} using Figure \ref{fig_CE}.

\bibliography{biblio}

\begin{thebibliography}{10}
\providecommand{\natexlab}[1]{#1}
\providecommand{\url}[1]{\texttt{#1}}
\expandafter\ifx\csname urlstyle\endcsname\relax
  \providecommand{\doi}[1]{doi: #1}\else
  \providecommand{\doi}{doi: \begingroup \urlstyle{rm}\Url}\fi

\bibitem[Colombo et~al.(2012)Colombo, Maathius, Kalisch, and
  Richardson]{Colombo12}
Diego Colombo, Marloes Maathius, Marcus Kalisch, and Thomas Richardson.
\newblock Learning high-dimensional directed acyclic graphs with latent and
  selection variables.
\newblock \emph{Annals of Statistics}, 40\penalty0 (1):\penalty0 294--321,
  April 2012.
\newblock \doi{10.1214/11-AOS940}.
\newblock URL \url{http://projecteuclid.org/euclid.aos/1333567191}.

\bibitem[Cowell et~al.(1999)Cowell, Lauritzen, David, and
  Spiegelhalter]{Cowell99}
Robert~G. Cowell, Steffen~L. Lauritzen, A.~Philip David, and David~J.
  Spiegelhalter.
\newblock \emph{Probabilistic Networks and Expert Systems}.
\newblock Springer-Verlag, Berlin, Heidelberg, 1st edition, 1999.
\newblock ISBN 0387987673.

\bibitem[Dagum et~al.(1992)Dagum, Galper, and Horvitz]{Dagum92}
Paul Dagum, Adam Galper, and Eric Horvitz.
\newblock Dynamic network models for forecasting.
\newblock In \emph{Proceedings of the Eighth International Conference on
  Uncertainty in Artificial Intelligence}, UAI'92, pages 41--48, San Francisco,
  CA, USA, 1992. Morgan Kaufmann Publishers Inc.
\newblock ISBN 1-55860-258-5.
\newblock URL \url{http://dl.acm.org/citation.cfm?id=2074540.2074546}.

\bibitem[Dagum et~al.(1995)Dagum, Galper, Horvitz, and Seiver]{Dagum95}
Paul Dagum, Adam Galper, Eric Horvitz, and Adam Seiver.
\newblock Uncertain reasoning and forecasting.
\newblock \emph{International Journal of Forecasting}, 11:\penalty0 73--87,
  1995.

\bibitem[Forr{\'e} and Mooij(2017)]{Forre17}
Patrick Forr{\'e} and Joris~M. Mooij.
\newblock Markov properties for graphical models with cycles and latent
  variables.
\newblock \emph{arXiv.org preprint}, arXiv:1710.08775 [math.ST], October 2017.
\newblock URL \url{https://arxiv.org/abs/1710.08775}.

\bibitem[Lauritzen et~al.(1990)Lauritzen, Dawid, Larsen, and
  Leimer]{Lauritzen90}
S.~L. Lauritzen, A.~P. Dawid, B.~N. Larsen, and H.~G. Leimer.
\newblock {Independence Properties of Directed Markov Fields}.
\newblock \emph{Networks}, 20\penalty0 (5):\penalty0 491--505, August 1990.
\newblock \doi{10.1002/net.3230200503}.
\newblock URL \url{http://dx.doi.org/10.1002/net.3230200503}.

\bibitem[Spirtes(1994)]{Spirtes94}
Peter Spirtes.
\newblock Conditional independence properties in directed cyclic graphical
  models for feedback.
\newblock Technical report, Carnegie Mellon University, 1994.

\bibitem[Spirtes(2001)]{Spirtes01}
Peter Spirtes.
\newblock An anytime algorithm for causal inference.
\newblock In \emph{in the Presence of Latent Variables and Selection Bias in
  Computation, Causation and Discovery}, pages 121--128. MIT Press, 2001.

\bibitem[Strobl(2018)]{Strobl18}
Eric~V. Strobl.
\newblock A constraint-based algorithm for causal discovery with cycles, latent
  variables and selection bias.
\newblock \emph{International Journal of Data Science and Analytics}, Nov 2018.
\newblock ISSN 2364-4168.
\newblock \doi{10.1007/s41060-018-0158-2}.
\newblock URL \url{https://doi.org/10.1007/s41060-018-0158-2}.

\bibitem[Strobl(2019)]{Strobl19}
Eric~V. Strobl.
\newblock Improved causal discovery from longitudinal data using a mixture of
  dags.
\newblock In Thuc~Duy Le, Jiuyong Li, Kun Zhang, Emre K{\i}c{\i}man~Peng Cui,
  and Aapo Hyv\"{a}rinen, editors, \emph{Proceedings of Machine Learning
  Research}, volume 104 of \emph{Proceedings of Machine Learning Research},
  pages 100--133, Anchorage, Alaska, USA, 05 Aug 2019. PMLR.
\newblock URL \url{http://proceedings.mlr.press/v104/strobl19a.html}.

\end{thebibliography}

\end{document}